\documentclass[11pt,leqno]{amsart}
\usepackage{latexsym,amssymb,amscd}

\def\demo{\noindent{\bf Proof. }}
\def\sqr#1#2{{\vcenter{\hrule height.#2pt
        \hbox{\vrule width.#2pt height#1pt \kern#1pt
                \vrule width.#2pt}
        \hrule height.#2pt}}}
\def\square{\mathchoice\sqr64\sqr64\sqr{4}3\sqr{3}3}
\def\QED{\hfill$\square$}

\newtheorem{lem}{Lemma}[section]
\newtheorem{prop}[lem]{Proposition}
\newtheorem{cor}[lem]{Corollary}
\newtheorem{thm}[lem]{Theorem}
\newtheorem{ques}[lem]{Question}
\newtheorem{remark}[lem]{Remark}
\newtheorem{dis}[lem]{Discussion}
\newtheorem{exam}[lem]{Example}

\DeclareMathOperator{\n}{\mathbf n}
\DeclareMathOperator{\m}{\mathbf m}

\DeclareMathOperator\depth{depth}
\DeclareMathOperator\Tor{Tor}
\DeclareMathOperator\Soc{Soc}

\DeclareMathOperator\Sym{ Sym}
\DeclareMathOperator\Hom{Hom}
\DeclareMathOperator\Ext{Ext}

\DeclareMathOperator\Quot{Quot}
\DeclareMathOperator{\hgt}{ht}
\DeclareMathOperator{\pd}{pd}

\DeclareMathOperator{\Z}{\mathbb Z}

\DeclareMathOperator{\F}{\mathbb F}

\DeclareMathOperator{\A}{\mathcal A}
\DeclareMathOperator{\B}{\mathcal B}

\DeclareMathOperator{\M}{\mathcal M}
\DeclareMathOperator{\R}{\mathcal R}

\def\alert#1{\smallskip{\hskip\parindent\vrule%
\vbox{\advance\hsize-2\parindent\hrule\smallskip\parindent.4\parindent%
\narrower\noindent#1\smallskip\hrule}\vrule\hfill}\smallskip}

\begin{document}

\title[]{The Gorenstein and complete intersection \\
properties of associated graded rings}

\date{\today }

\author{\ }

\address{\ }

\thanks{ Mee-Kyoung Kim is
partially supported by the Korea Science and Engineering
Foundation (R04-2003-000-10113-0). Bernd Ulrich is partially
supported by the National Science Foundation (DMS-0200858). }

\dedicatory{Dedicated to Wolmer Vasconcelos on the occasion of his
65th birthday.}


\maketitle

\begin{center}
\begin{tabular}{c}
{\bf William  Heinzer} \quad and \quad {\bf Bernd Ulrich} \\
Department of Mathematics, Purdue University \\
West Lafayette, Indiana 47907, USA\\
E-mail: heinzer@math.purdue.edu \\
E-mail: ulrich@math.purdue.edu \\
\ \\
{\bf Mee-Kyoung Kim} \\
Department of Mathematics, Sungkyunkwan University \\
Jangangu Suwon 440-746, Korea \\
E-mail: mkkim@math.skku.ac.kr \\
\end{tabular}
\end{center}

\bigskip

\begin{abstract}
\noindent
Let $I$ be an $\m$-primary ideal of a
Noetherian local ring $(R,\m)$.  We
consider the Gorenstein and complete
intersection properties   of the associated graded
ring $G(I)$ and the fiber cone $F(I)$ of $I$ as reflected
in their defining ideals as homomorphic images of
polynomial rings over $R/I$ and $R/\m$ respectively.
In case  all the higher conormal modules of $I$ are free over
$R/I$, we observe that: (i) $G(I)$ is Cohen-Macaulay iff $F(I)$ is
Cohen-Macaulay, (ii) $G(I)$ is Gorenstein iff both $F(I)$ and
$R/I$ are Gorenstein, and (iii) $G(I)$ is a relative complete
intersection iff  $F(I)$ is a complete intersection.
In case  $(R,\m)$ is Gorenstein,
we give a necessary and sufficient condition for $G(I)$ to
be Gorenstein in terms of residuation of powers of $I$
with respect to a reduction $J$ of $I$ with $\mu(J) = \dim R$ and the
reduction number $r$  of $I$ with respect to $J$.  We
prove that $G(I)$ is Gorenstein if and only if
$J:I^{r-i} = J + I^{i+1}$  for $0 \leq i \leq r-1$.
If $(R,\m)$ is a Gorenstein local ring and $I \subseteq \m$
is an ideal having a reduction $J$ with reduction
number $r$ such that
$\mu(J) = \hgt(I) = g > 0$, we prove that
the extended Rees algebra $R[It, t^{-1}]$ is
quasi-Gorenstein with {\bf a}-invariant $a$ if and only if
$J^i:I^r = I^{i+a-r+g-1}$ for every $i \in \Z$. If,
in addition, $\dim R = 1$, we show that $G(I)$ is
Gorenstein if and only if $J^i:I^r = I^i$
for $1 \leq i \leq r$.
\end{abstract}

\baselineskip=15pt

\section{Introduction}

For an ideal $I \subseteq \m$ of a Noetherian local ring $(R,\m)$,
several graded rings naturally associated to $I$ are:
\begin{enumerate}
\item
the symmetric algebra $\Sym_R(I)$ and
the Rees algebra $\R = R[It] = \bigoplus_{i \ge 0} I^it^i$ \, (considered as
a subalgebra of the poynomial ring $R[t]$),
\smallskip
\item
the extended symmetric algebra $\Sym_R(I, t^{-1})$ and
the extended Rees algebra $R[It, t^{-1}] = \bigoplus_{i \in \Z}I^it^i$
\, (using the convention that $I^i=R$ for $i \leq 0$),
\smallskip
\item
the symmetric algebra $\Sym_{R/I}(I/I^2)$ and
the associated graded ring $G(I) = R[It, t^{-1}]/(t^{-1}) = R[It]/IR[It] =
\bigoplus_{i \ge 0}I^i/I^{i+1}$, and
\smallskip
\item
the  symmetric algebra $\Sym_{R/\m}(I/\m I)$ and the fiber cone
$F(I) = \newline R[It]/\m R[It] = \bigoplus_{i \ge 0}I^i/\m I^i$.
\end{enumerate}

These graded rings encode information about $I$ and its powers.
The {\em analytic spread } of $I$, denoted $\ell(I)$, is the
dimension of the fiber cone $F(I)$ of $I$. An ideal $J \subseteq
I$ is said to be a {\em reduction } of $I$ if there exists a
nonnegative integer $k$ such that $JI^k = I^{k+1}$. It then
follows that $J^jI^k = I^{k+j}$ for every nonnegative integer $j$.
These concepts were introduced by Northcott and Rees in \cite{NR}.
If $J$ is a reduction of $I$, then $J$  requires at least
$\ell(I)$ generators. Reductions of $I$ with $\ell(I)$
generators are necessarily minimal reductions in the sense that no
properly smaller ideal is a reduction of $I$. They correspond to
Noether normalizations of $F(I)$ in the sense that $a_1, \dots,
a_{\ell}$ generate a reduction of $I$ with $\ell = \ell(I)$
if and only if their images
$\overline{a_i} \in I/\m I \subseteq F(I)$ are algebraically
independent over $R/\m$ and $F(I)$ is integral over the polynomial
ring $(R/\m)[\overline{a_1}, \dots , \overline{a_{\ell}}]$. In
particular, if $R/\m$ is infinite, then there exist
reductions of $I$ generated by $\ell(I)$ elements.

Suppose $(R,\m)$ is a Gorenstein local ring of dimension $d$ and
$I$ is an $\m$-primary ideal. We are interested in conditions for
the associated graded ring $G(I)$ or the fiber cone $F(I)$ to be
Gorenstein. Assume $J = (a_1, \ldots, a_d)R$ \, is a reduction of
$I$. If $G(I)$ is Cohen-Macaulay, then the images $a_1^*, \ldots,
a_d^*$ of $a_1, \ldots, a_d$ in $I/I^2$ form a regular sequence on
$G(I)$, and $G(I)$ is Gorenstein if and only if $G(I)/(a_1^*,
\ldots, a_d^*)G(I)$ is Gorenstein. Write $\overline{R} = R/(a_1,
\ldots, a_d)R$ and $\overline{I} =I\overline{R}$. Then
$\overline{R}$ is a zero-dimensional Gorenstein local ring and
$G_{\overline{R}}(\overline{I}) = G(I)/(a_1^*, \ldots,
a_d^*)G(I)$. Thus, under the hypothesis that $G(I)$ is
Cohen-Macaulay,
 the question of whether
$G(I)$ is Gorenstein reduces to a zero-dimensional setting.

For $J \subseteq I$ a reduction of $I$,
the {\em reduction number } $r_J(I)$ of $I$ with
respect to $J$ is the smallest nonnegative integer $r$
such that $JI^r = I^{r+1}$. If $I$ is an $\m$-primary
ideal of  a $d$-dimensional Gorenstein local ring $(R,\m)$
and $J$ is a $d$-generated reduction of $I$,
then the reduction number $r = r_J(I)$ plays an
important role in considering the Gorenstein property
of $G(I)$. If $r = 0$, then $J = I$ is generated by
a regular sequence and $G(I)$ is a polynomial ring in $d$
variables over the zero-dimensional Gorenstein local ring
$R/I$. Thus $G(I)$ is Gorenstein in this case.
If $r = 1$ and $d \ge 2$, then a result of Corso-Polini
\cite[Cor.3.2]{CP} states that $G(I)$ is Gorenstein
if and only if
$J:I = I$, that is if and only if the ideal $I$ is self-linked
with respect to the minimal reduction parameter ideal $J$.

In Theorem \ref{3.8} we extend this result of Corso-Polini on the
Gorenstein property of $G(I)$ in case $I$ is an $\m$-primary
ideal. We prove that $G(I)$ is Gorenstein if and only if
$J:I^{r-i} = J + I^{i+1}$ for $0 \leq i \leq r-1$. This also gives
an analogue for Gorenstein rings  to a well-known result about
Cohen-Macaulay rings (see  \cite[Lemma 2.2]{HM})  that asserts:
Suppose $I$ is an $\m$-primary ideal of a Cohen-Macaulay local
ring $(R,\m)$ of dimension $d$ and the elements $a_1, \ldots, a_s$
are a  superficial sequence for $I$. Let $\overline{R} = R/(a_1,
\ldots, a_s)R$ \, and $\overline{I} = I\overline{R}$. If $s \le
d-1$, then $G(\overline{I})$ is Cohen-Macaulay implies $G(I)$ is
Cohen-Macaulay. Corollary~\ref{3.9} says that even for $s = d$ it
holds that $G(\overline{I})$ is Gorenstein implies $G(I)$ is
Gorenstein if one assumes that $I^r \not\subseteq J = (a_1,
\ldots, a_d)R$.

Let  $B = \oplus_{i \in \Z}B_i$ be a Noetherian $\Z$-graded ring which
is *local in the sense that it has a unique maximal homogeneous
ideal $\M$ \cite[(1.5.13), p.35]{BH}. Notice that $R:=B_0$
is a Noetherian local ring and $B$ is finitely generated as an $R$-algebra
\cite[Thm.1.5.5, p.29]{BH}. We assume for simplicity that $R$ is Gorenstein. Consider
a homogeneous presentation $B \cong S/H$ with $S=R[X_1, \ldots, X_n]$
a $\Z$-graded polynomial ring and
$H$ a homogeneous ideal of height $g$. Let $\sigma \in \Z$ be the
sum of the degrees of the variables $X_1, \ldots, X_n$. We write
$\omega_B= \Ext_S^g(B,S)(-\sigma)$ and call this module the
{\it graded canonical module} of $B$. One easily sees that $\omega_B$
is a finitely generated graded $B$-module that is uniquely
determined up to homogeneous $B$-isomorphisms.
The ring $B$ is said to be {\it quasi-Gorenstein} in case $\omega_B \cong B(a)$
for some $a \in \Z$. If the maximal homogeneous ideal $\M$ of $B$
is a maximal ideal, then the integer $a$ is
uniquely determined and is called the {\it {\bf a}-invariant} of
$B$. We will use the following facts, which are
readily deduced from the above definition of graded canonical modules:

$\bullet$ The localization $(\omega_B)_{\M}$ is the canonical module of
the local ring $B_{\M}$.

$\bullet$ The module $\omega_B$ satisfies $S_2$.

$\bullet$ The ring $B$ is (locally) Gorenstein if and only if it
is quasi-Gorenstein and (locally) Cohen-Macaulay.

$\bullet$ Let $A$ be a $\Z$-graded subring of $B$ with unique maximal
homogeneous ideal $\M \cap A$ and $A_0=R=B_0$, so that $A$ is
Cohen-Macaulay and $B$ is a
finitely generated $A$-module; then $\omega_B \cong \Hom_A(B,\omega_A)$
as graded $B$-modules.

Notice that a quasi-Gorenstein ring
necessarily satisfies $S_2$. Thus if $B$ is quasi-Gorenstein and
$\dim B \le 2$, then $B$ is Gorenstein. In higher dimensions there
do exist examples of quasi-Gorenstein rings that are not
Gorenstein. There exists an  example of a prime ideal $P$ of
height two in a $5$-dimensional regular local ring $R$ such that
the extended Rees algebra $R[Pt, t^{-1}]$ is  quasi-Gorenstein but
not Gorenstein \cite[Ex.4.7]{HH}. We are interested in
classifying quasi-Gorenstein extended Rees algebras  and saying
more about when such rings are Gorenstein. In Theorem~\ref{4.1},
we prove that if $(R,\m)$ is a Gorenstein local ring and $I
\subseteq \m$ is an ideal having a reduction $J$ with reduction
number $r$ such that $\mu(J) = \hgt(I) = g > 0$, then the extended
Rees algebra $R[It, t^{-1}]$ is quasi-Gorenstein with {\bf
a}-invariant $a$ if and only if $J^i:I^r = I^{i+a-r+g-1}$ for
every $i \in \Z$. A natural question here that we consider but
resolve only in special cases is whether $R[It, t^{-1}]$ is
Gorenstein if it is  quasi-Gorenstein.  Equivalently, is the
associated graded ring $G(I)$ then Gorenstein. We observe that
this is true if $\dim R = 1$ or if $\dim R = 2$ and $R$ is
regular. If  $\dim R = 1$, Corollary \ref{4.4} implies that $G(I)$
is Gorenstein if and only if $J^i:I^r = I^i$ for $1 \leq i \leq
r$.

\medskip

\section{Defining ideals and freeness of the higher
conormal modules}

Let $(R, \m)$ be a Noetherian local ring and let
$I = (a_1, \ldots, a_n)R \subseteq \m$
be an ideal of $R$. Consider the  presentation of the
Rees algebra $R[It]$ as a homomorphic image of a polynomial
ring over $R$ obtained  by defining an $R$-algebra homomorphism
$\tau : R[X_1, \ldots , X_n] \to R[It]$ such that
$\tau(X_i) = a_it$ for $1 \leq i \leq  n$. Now
define $\psi = \tau \otimes R/I : (R/I)[X_1, \ldots, X_n] \to G(I)$,
where $\psi(X_i) = a_i + I^2 \in I/I^2 = G_1$, and
$\phi = \tau \otimes R/\m : (R/\m)[X_1, \ldots, X_n] \to F(I)$,
where $\phi(X_i) = a_i + \m I \in  I/\m I = F_1$.

We have  the following commutative diagram for which the rows are
exact and the column maps are surjective.

\begin{equation}\label{eq1}\tag{2.1}
\text{$\CD
0 @>>> E := \ker(\tau)  @>>> R[X_1, \ldots, X_n] @>\tau>> R[It] @>>> 0 \\
@. @V\pi_1 VV        @V\pi_2 VV                @V\pi_3 VV \\
0 @>>> L := \ker(\psi)  @>>> (R/I)[X_1, \ldots, X_n] @>\psi>> G(I)
@>>> 0 \\
@. @V\pi_1' VV        @V\pi_2' VV                @V\pi_3' VV \\
0 @>>> K := \ker(\phi) @>>> (R/\m)[X_1, \ldots, X_n] @>\phi>> F(I) @>>> 0 .
\endCD $}
\end{equation}

\vskip 3ex

The ideal $E = 0$ if and only if $n = 1$
and $a_1$ is a regular element of $R$, while
$L  = 0$ if and only if $a_1, \ldots, a_n$
form a regular sequence, see for example
\cite[Cor.5.13, p.154]{K}. A necessary and
sufficient condition for $K = 0$ is
that $a_1, \ldots, a_n$ be what Northcott and Rees
\cite{NR} term {\it analytically independent. }

Let $v$ be the maximal degree of a homogeneous minimal generator
of the ideal $E$. The integer $N_{\R}(I):=\max\{1,v\}$
is called the {\it relation type } of the Rees algebra $\R =
R[It]$ with respect to the given generating set $a_1, \ldots,
a_n$. The relation type of $R[It]$ may also be defined by
considering  the kernel $N$ of the canonical homomorphism from the
symmetric algebra $\Sym_R(I)$ onto $R[It]$; then  $N_{\R}(I) =
\max\{1, w \}$, where $w$ denotes the maximal degree of a
homogeneous minimal generator of $N$.  Since the symmetric algebra
$\Sym_R(I)$ is independent of the choice of generators for $I$, it
follows that the relation type  $N_{\R}(I)$ of $R[It]$ is
independent of the generating set of $I$. The ideal $I$ is said to
be of {\it linear type} if $N_{\R}(I) = 1$. Thus $I$ is of linear
type if and only if $R[It]$ is canonically isomorphic to the
symmetric algebra $\Sym_R(I)$ of $I$.

The relation type $N_G(I)$ of the associated graded ring $G(I)$
is  defined in a similar
way, using the ideal $L$ or the kernel of the
canonical homogeneous epimorphism
$\beta : \Sym_{R/I}(I/I^2)  \to G(I)$.
Likewise, the relation type $N_F(I)$ of the fiber cone $F(I)$
is defined via the ideal $K$ or the kernel of the
canonical homogeneous
epimorphism $\gamma : \Sym_{R/\m}(I/\m I) \to F(I)$.
The surjectivity of the maps $\pi_1$ and $\pi_1'$ in diagram
(\ref{eq1}) imply that the inequalities $N_F(I) \le N_G(I) \le
N_{\R}(I)$ hold in general.

\setcounter{lem}{1}

\begin{dis}\label{2.1} {\em
It is shown by Valla \cite[Thm.1.3]{Va} and by
Herzog-Simis-Vasconcelos \cite[Thm.3.1]{HSV} that if $N_G(I) = 1$,
then also $N_{\R}(I) = 1$, that is, the relation type of $G(I)$ is
one  if and only if $I$  is of linear type. Planas-Vilanova shows
\cite[Prop.3.3]{PV} that the equality $N_{\R}(I) = N_G(I)$ also
holds in the case where $N_G(I) > 1$. We reprove this fact by
using an `extended symmetric algebra' analogue to the extended
Rees algebra $R[It, t^{-1}]$ of the ideal $I$. As an $R$-module we
define the {\it extended symmetric algebra} \, $\Sym_R(I, t^{-1})$
of $I$ to be $\oplus_{i \in \Z}S_i(I)$, where $S_i(I) = t^iR \cong
R$ for $i \le 0$ and $S_i(I)$ is the $i$-th symmetric power of $I$
for $i \ge 1$. We define a multiplication on $\Sym_R(I, t^{-1})$
that extends the ring structure on the standard symmetric algebra
$\Sym_R(I) = \oplus_{i \ge 0}S_i(I)$. To do this, it suffices to
define multiplication by $t^{-1}$ on $S_m(I)$ for $m \ge 1$.
Consider the product $I \times \cdots \times I$ of $m$ copies of
the ideal $I$ and the map $I \times \cdots \times I \to
S_{m-1}(I)$ that takes $(a_1, \ldots, a_m) \mapsto (a_1\cdot
a_2\cdot \ldots \cdot a_{m-1})a_m$, where the multiplication on
$a_1, \ldots ,a_{m-1}$ is multiplication in the symmetric algebra
and where the multiplication with $a_m$ is the scalar
multiplication given by the $R$-module structure of $S_{m-1}(I)$.
This map is $m$-multilinear over $R$ and hence factors through the
tensor power $I \otimes \ldots \otimes I$. Moreover, the map is
symmetric in $a_1, \ldots, a_m$ and hence  induces a map $S_m(I)
\to S_{m-1}(I)$. To verify this symmetry, it suffices to observe
that by associativity and commutativity of the two
multiplications, one has $(a_1\cdot  \ldots \cdot a_{m-2}\cdot
a_{m-1})a_m  = (a_1\cdot  \ldots \cdot a_{m-2})\cdot (a_{m-1}a_m)
= (a_1\cdot  \ldots \cdot a_{m-2})\cdot (a_ma_{m-1})  = (a_1\cdot
\ldots \cdot a_{m-2}\cdot a_m)a_{m-1}$. The multiplication by
$t^{-1}$ on $S_m(I)$ we just defined coincides with the
downgrading homomorphism $\lambda_{m-1}$ introduced by Herzog,
Simis and Vasconcelos \cite[p.471]{HSV}. Assigning to the element
$t^{-1}$ degree $-1$, $\Sym_R(I, t^{-1})$ becomes a $\Z$-graded
$R$-algebra. This algebra is $^*$local \cite[(1.5.13), p.35]{BH}
in the sense that the ideal of $\Sym_R(I, t^{-1})$ generated by
$t^{-1}, \m$ and $\oplus_{i \ge 1}S_i(I)$ is the unique maximal
homogeneous ideal of $\Sym_R(I, t^{-1})$.

The canonical surjective $R$-algebra homomorphism from the
symmetric algebra $\Sym_R(I)$ onto the Rees algebra $R[It]$
extends to a surjective homogeneous $R$-algebra homomorphism
$\alpha : S = \Sym_R(I, t^{-1}) \to R[It, t^{-1}]$, where
$\alpha(t^{-1}) = t^{-1}$. Notice that the two maps have the same
kernel $\A$. Tensoring the short exact sequence
\[
\CD
0 @>>> \A  @>>> S
@>\alpha>> R[It, t^{-1}] @>>> 0 \\
\endCD
\]
with $\otimes_SS/t^{-1}S$ gives the following isomorphisms:
\[
\CD
0 @>>> \A \otimes S/t^{-1}S  @>>> S \otimes S/t^{-1}S
@>>> R[It, t^{-1}]\otimes S/t^{-1}S @>>> 0 \\
@. @V\cong VV        @V\cong VV                @V\cong VV \\
0 @>>> \B: = \ker(\beta)@>>> \Sym_{R/I}(I/I^2) @>\beta>> G(I)  @>>> 0 .\\
\endCD
\]
Here we are using that $t^{-1}$ is a regular element on $R[It,
t^{-1}]$. A graded version of Nakayama's lemma \cite[(1.5.24),
p.39]{BH} now implies that $\A$ is generated in degrees $\leq
N_G(I)$ as a module over $S=\Sym_R(I)[t^{-1}]$ and hence over
$\Sym_R(I)$. Therefore the equality $N_{\R}(I) = N_G(I)$ holds in
general.}
\end{dis}

\medskip

We observe in Corollary \ref{2.6}
a sufficient condition for $N_G(I) = N_F(I)$.
There exist, however,   examples where
$N_G(I) > N_F(I)$. If $(R,\m)$ is a $d$-dimensional Noetherian
local ring and $I = (a_1, \ldots, a_d)R$ is an $\m$-primary ideal,
then $N_F(I) = 1 $ \cite[Thm.14.5,
p.107]{M}; while if $R$ is not Cohen-Macaulay  it may happen
that $N_G(I) > 1$. For example, let $k$ be a field and let
$A = k[x,y] = k[X,Y]/(X^2Y,Y^3)$. Then $R = A_{(x,y)A}$
is a one-dimensional Noetherian local ring
and $I = xR$ is primary to the maximal ideal $\m = (x,y)R$. Let
$U$ be an indeterminate over $R/I$ and consider
the graded homomorphism $\psi: (R/I)[U] \to G(I)$,
where $\psi(U) = x + I^2 \in I/I^2 = G_1$. Let
$L = \oplus_{i \ge 0}L_i = \ker(\psi)$.
Since $R/I \cong k[X, Y]/(X, Y^3) \cong k[Y]/(Y^3)$,
we see that  $\lambda(R/I) = 3$, where $\lambda$ denotes
length. Similarly,
$R/I^2 \cong k[X, Y]/(X^2, Y^3)$ and
$\lambda(R/I^2) = 6$. Thus  $\lambda(I/I^2) = 3$. It
follows that $L_0 = L_1  = 0$. On the other hand,
$0 \neq (y+I)U^2 \in L_2$. Therefore $N_G(I) \ge 2$, while
$N_F(I) = 1$.

In the case where $N_F(I) \ge 2$, it would be
interesting to have necessary and sufficient conditions
for $N_G(I) = N_F(I)$.
That this is not always true is
shown, for example, by taking $k$ to be a field and
$A = k[x, y, z]=k[X,Y,Z]/(X^2,Y^2,XYZ^2)$.
Let $R = A_{(x,y,z)A}$ and let
$I = (y, z)R$. Then $F(I) \cong k[Y, Z]/(Y^2)$
has relation type 2, while
if $\psi : (R/I)[U, V]  \to G(I)$ with $\psi(U) = y +I^2$
and $\psi(V) = z + I^2$ is a presentation
of $G(I)$, then the relation $xyz^2 = 0$ in the
definition of $R$ implies $(x+I)UV^2 \in \ker(\psi)$.
Let $L = \oplus_{i \geq 0} L_i = \ker(\psi)$.
Since $R/I \cong k[X]/(X^2)$, we see that
$\lambda(R/I) = 2$.  We have
$\lambda(I/I^2) = 4$, hence  $L_1 = 0$. Also
$\lambda(I^2/I^3) = 4$, so $\lambda(L_2) = 2$
and $L_2$ is generated by $U^2$. Since
$(x+I)UV^2 \not\in (U^2+I)(R/I)[U, V]$, we see that
$N_G(I) \ge 3$.

\begin{thm}\label{2.3}
Let $(R,\m)$ be a Noetherian local ring and let $I$ be an
$\m$-primary ideal with $\mu(I) = \lambda(I/\m I) = n$. Let  $A =
(R/I)[X_1, \ldots, X_n]$ and let $\psi : A \to G(I)$ be defined as
in diagram $(\ref{eq1})$. The following are equivalent $:$
\begin{enumerate}
\item
$I^i/I^{i+1}$ is free as an $(R/I)$-module for every $i \geq 0$.
\item
$G(I)$ has finite projective dimension as an $A$-module,
i.e., $\pd_AG(I) < \infty$.
\item If $\F_\bullet$ is a homogeneous minimal free resolution of $G(I)$
as an $A$-module and $k = R/\m$,
then $\F_\bullet \otimes \, k$ is a homogeneous minimal free
resolution of $F(I) = G(I) \otimes  k$
as a module over $A \otimes k = k[X_1, \ldots, X_n]$.
\end{enumerate}
If these equivalent conditions hold, then $\depth G(I) = \depth F(I)$.
\end{thm}

\demo
$(1) \implies (3)$: Condition (1) implies that $\Tor^{R/I}_i(G(I), k) = 0$ for
every $i > 0$, and this implies condition (3).

\noindent
$(3) \implies (2)$: Let $\F_\bullet$ be a homogeneous minimal free resolution of
$G(I)$ over $A$. Condition (3) implies that $\F_\bullet \otimes \, k$ is
a homogeneous minimal free resolution of $F(I)$ over the regular ring $B = k[X_1, \ldots, X_n]$.
Since $B$ is regular, this homogeneous minimal free resolution of $F(I)$ is finite.
Therefore $\F_\bullet$ is finite and $\pd_AG(I) < \infty$.

\noindent $(2) \implies (1)$: Condition (2) implies that
$\pd_{R/I}(I^i/I^{i+1}) < \infty$ for every $i \geq 0$. Since
$\dim(R/I) = 0$, $I^i/I^{i+1}$ is free over $R/I$ by the
Auslander-Buchsbaum formula, see for instance \cite[Thm.1.3.3,
p.17]{BH}.

If these equivalent conditions hold, then we have
$\pd_AG(I) = \pd_BF(I)$ by (3). Hence
$\depth G(I) = \depth F(I)$ by the Auslander-Buchsbaum formula.
\QED

\begin{remark}\label{2.5}{\em
If $(R, \m)$ is a one-dimensional Cohen-Macaulay local ring, $I$
is an $\m$-primary ideal  and $J=xR$ is a principal reduction of
$I$ with reduction number $r = r_J(I)$, then $x^{i-r}I^r = I^i$
for $i \ge r$, so $I^i/I^{i+1} \cong I^r/I^{r+1}$. Thus in this
case a fourth equivalent condition  in Theorem \ref{2.3} is that
$I^i/I^{i+1}$ is free over $R/I$ for every $i \le r$. On the other
hand, if $I \subset \m$ is an ideal of a zero-dimensional
Noetherian local ring $(R,\m)$, then a condition equivalent  to
the conditions of Theorem \ref{2.3} is that $\lambda(R) =
\lambda(R/I)\cdot \lambda(F(I))$. }
\end{remark}

Under the equivalent conditions of Theorem \ref{2.3}, the
ring $R$ is said to be {\it normally flat} along the ideal $I$.
This is a concept introduced by Hironaka in his work on resolution
of singularities \cite[p.188]{M}. A well known
result of Ferrand \cite{F} and Vasconcelos \cite[Cor.1]{V} asserts that
if the conormal module $I/I^2$ is a free module over $R/I$
and if $R/I$ has finite
projective dimension as an $R$-module, then $I$ is generated
by a regular sequence. It then  follows that $N_G(I) = 1$. Thus,
in the case where $R$ is a regular local ring, the equivalent
conditions of Theorem \ref{2.3} imply $\mu(I) = \dim R$.
However, as we indicate in Example \ref{2.9},  there exist  examples of
Gorenstein local rings $(R,\m)$ having
$\m$-primary ideals $I$  such that
$\mu(I) > \dim R$ and the equivalent
conditions of Theorem \ref{2.3} hold.

In Examples \ref{2.9}, \ref{2.13}, \ref{2.10}, \ref{3.5},
\ref{3.11}, \ref{4.5}, and \ref{4.51},
we present examples
involving an additive monoid $S$ of the nonnegative
integers that contains all sufficiently large integers
and a complete one-dimensional local
domain of the form $R = k[[t^s \, | \,  s \in S]]$.
The formal power series ring $k[[t]] = R[t]$ is the
integral closure of $R$ and is a finitely
generated $R$-module. Properties of $R$ are closely
related to properties of the numerical semigroup $S$.
For example, $R$ is Gorenstein if and only if $S$ is
symmetric \cite[Thm.4.4.8, p.178]{BH}.

\begin{exam}\label{2.9}{\em
Let $k$ be a field, $R = k[[t^4, t^9, t^{10}]]$ and $I =
(t^8, t^9, t^{10})R$. Then $R$ is a one-dimensional
Gorenstein local domain, $\mu(I) = 3$ and
$J = t^8R$ is a reduction of $I$.
We have $JI \subsetneq I^2$ and $JI^2 = I^3$.
Hence $r_J(I) = 2$. If  $w$ denotes the image of
$t^4$ in $R/I$, then  $R/I \cong k[w]$, where
$w^2 = 0$. Thus $R/I$ is Gorenstein with
$\lambda(R/I) = 2$.
Let $A = (R/I)[X, Y, Z]$ and define
$\psi : A \to G(I)$ by $\psi(X) = t^8 + I^2,
\psi(Y) = t^9 + I^2$ and $\psi(Z) = t^{10} + I^2$.
Consider the  short exact sequence
$$
\CD
0 @>>> L = \ker(\psi)  @>>> A
@>\psi>> G(I) @>>> 0.
\endCD
$$
The ring $G(I) \cong A/L$ has multiplicity $e(I) = 8$, $A$
is a Cohen-Macaulay ring of multiplicity 2,
and the two quadrics
$XZ - Y^2, wX^2 - Z^2$ form an $A$-regular
sequence contained in $L$.
Hence $L = (XZ - Y^2, wX^2 - Z^2)A$.
Therefore  the equivalent conditions of
Theorem \ref{2.3} are satisfied.  It also follows that
$I^i/I^{i+1}$ is free of rank 4 over $R/I$ for every
$i \ge 2$.}
\end{exam}

\begin{cor}\label{2.6}
Let $(R,\m)$ be a Noetherian  local ring and  let $I$ be an
$\m$-primary ideal with $\mu(I) = \lambda(I/\m I) = n$. With
notation as in Theorem $\ref{2.3}$ and diagram $(\ref{eq1})$, let
$A=(R/I)[X_1, \ldots, X_n]$ and $B=(R/\m)[X_1, \ldots, X_n]$, and
let $\M_G = (\m/I, G(I)_+)$ and $\M_F = F(I)_+$ denote the maximal
homogeneous ideals of $G(I)$ and $F(I)$, respectively. Suppose the
equivalent conditions of Theorem $\ref{2.3}$ hold. Then$:$
\begin{enumerate}
\item
$G(I)$ is Cohen-Macaulay $\Longleftrightarrow$  $F(I)$ is Cohen-Macaulay.

\item
If $G(I)$, or equivalently $F(I)$, is Cohen-Macaulay,
then the type of $G(I)_{\M_G}$ is the type of $R/I$ times the type
of $F(I)_{\M_F}$. In particular,
$G(I)$ is Gorenstein $\Longleftrightarrow$  both $F(I)$ and
$R/I$ are  Gorenstein.

\item
The relation type $N_G(I)$ of the associated graded ring $G(I)$
is equal to the
relation type $N_F(I)$ of the fiber cone $F(I)$.

\item
The defining ideal $L$ of $G(I)$ is generated by a regular sequence
on $A$ if and only if the defining ideal $K$ of $F(I)$ is generated
by a regular sequence on $B$.

\item
The multiplicity of $G(I)$ is $\lambda(R/I)\cdot e(F(I))$,
where $e(F(I))$ denotes the multiplicity of $F(I)$.
\end{enumerate}
\end{cor}

\demo
As $G(I)$ is flat over $R/I$, statements (1) and (2) follow from
\cite[Prop.1.2.16, p.13]{BH}.
Since the relation types of $F(I)$ and $G(I)$ are determined
by the degrees of minimal generators of  first syzygies,
statement (3) follows from part (3) of
Theorem \ref{2.3}.  Indeed, with the notation of diagram (\ref{eq1}),
the relation type of $G(I)$ is $\max\{1, w\}$, where $w$ is the
maximal degree of a homogeneous minimal generator of $\ker(\psi)$.
Since $\F_{\bullet} \otimes k$ is a minimal free resolution of $F(I)$,
$\pi_1'$ maps a set of homogeneous
minimal generators of $L = \ker(\psi)$ onto a set of homogeneous
minimal generators of $K = \ker(\phi)$. Therefore $N_G(I) = N_F(I)$.
Since $A$ and $B$ are Cohen-Macaulay and $\hgt(L) = \hgt(K)$,
it also follows that
$L$ is generated by a regular sequence on $A$  if and only if
$K$ is generated by a regular sequence on $B$, which is part (4). Statement (5) is
clear in view of the freeness of the $I^i/I^{i+1}$ over $R/I$. \QED

\medskip

We observe
in part (1) of Remark \ref{2.4} that $L$ is generated by
a regular sequence implies $K$ is generated by a regular sequence
holds even without the equivalent conditions of Theorem \ref{2.3}.

In Proposition \ref{2.5.5}, we give
a partial converse to part (5) of
Corollary \ref{2.6}. Proposition \ref{2.5.5} is closely
related to results of Shah \cite[Lemma 8 and Thm.8]{Sh}.

\begin{prop} \label{2.5.5} Let $(R,\m)$ be a Noetherian local
ring and let $I$ be an $\m$-primary ideal. With notation
as in Theorem $\ref{2.3}$ assume that
$e(G(I)) = \lambda(R/I) \cdot e(F(I))$.
\begin{enumerate}
\item
If all associated primes of $F(I)$ have the same dimension,
then $I^i/I^{i+1}$ is a free $(R/I)$-module for every $i \ge 0$.
\item
If all relevant associated primes of $F(I)$ have the same dimension,
then $I^i/I^{i+1}$ is a free $(R/I)$-module for all sufficiently
large $i$.
\end{enumerate}
\end{prop}

\demo
A composition series of $R/I$ induces a filtration on $G(I)$
whose factors are homogeneous $F(I)$-modules of the form
$F(I)/{\bf a}_j$ with $1 \le j \le \lambda(R/I)$. Since
$e(G(I)) = \lambda(R/I) \cdot e(F(I))$, these factors all
have the same dimension and the same multiplicity as $F(I)$.
Thus in the setting of part (1), ${\bf a}_j = 0$ for every $j$,
whereas in part (2), ${\bf a}_j$ is concentrated in finitely
many degrees for every $j$. Computing the Hilbert function of
$G(I)$ one sees that $\lambda(I^i/I^{i+1}) = \lambda(R/I) \cdot
\lambda(I^i/\m I^i) = \lambda(R/I) \cdot \mu(I^i)$ for every
$i \ge 0$  in the
setting of part (1), and for every $i >> 0$ in the setting of
part (2). \QED

\medskip

In Example \ref{2.13}, we exhibit an $\m$-primary ideal $I$ of a
one-dimensional Cohen-Macaulay local domain $(R,\m)$ such
that $I^i/I^{i+1}$ is free over $R/I$ for every $i \ge 2$,
while  $I/I^2$ is not free over $R/I$. This example illustrates
that the equality $e(G(I)) = \lambda(R/I) \cdot e(F(I))$
does not imply the equivalent conditions of Theorem \ref{2.3},
even if the assumption of part (2) of Proposition \ref{2.5.5} is
satisfied.

\begin{exam}\label{2.13} {\em
Let $k$ be a field, $R = k[[t^3, t^7, t^{11}]]$ and $I =
(t^6, t^7, t^{11})R$ as in \cite[Ex.6.4]{DV}.
Then  $\mu(I) = 3$ and
$J = t^6R$ is a principal reduction of $I$.
Also  $JI \subsetneq I^2$ and $JI^2 = I^3$, so
$r_J(I) = 2$.  We have
$\lambda(R/I) = 2$ and  $\lambda(R/I^2) = 7$.
Hence $\lambda(I/I^2) = 5$ and $I/I^2$ is not
free over $R/I$.  On the other hand,
$I^2 = (t^{12}, t^{13}, t^{14})R$. Therefore   $I^i/I^{i+1}$ is
generated by 3 elements and $\lambda(I^i/I^{i+1}) = 6$,
so $I^i/I^{i+1}$ is free over $R/I$ of rank 3 for
every $i \ge 2$.}
\end{exam}

By Proposition \ref{2.5.5},  the fiber cone $F(I)$
in Example \ref{2.13} is not Cohen-Macaulay. To
see this explicitly, let  $w$ denote the image of $t^3$ in $R/I$ and let
$A = (R/I)[X, Y, Z]$. Define
$\psi : A \to G(I)$ by $\psi(X) = t^6 + I^2,
\psi(Y) = t^7 + I^2$ and $\psi(Z) = t^{11} + I^2$,  and
consider the  exact sequence
$$
\CD
0 @>>> L = \ker(\psi)  @>>> A
@>\psi>> G(I) @>>> 0.
\endCD
$$
Then  $L = (wZ, XZ - wY^2, YZ, Z^2, Y^3 - wX^3)A$.
Hence
$$
e(G(I))=e(I)=6=\lambda(G(I)/\psi(X)G(I)),
$$
showing that $G(I)$ is
Cohen-Macaulay, see for instance \cite[Thm.17.11, p.138]{M}.
However, the fiber cone $F(I)$ is isomorphic
to $(R/\m)[X, Y, Z]/(XZ, YZ, Z^2, Y^3)$ and is
not Cohen-Macaulay.

\vskip 3ex
\begin{remark}\label{2.4}{\em
With notation as in Theorem \ref{2.3} and diagram (\ref{eq1}), we have:
\begin{enumerate}
\item
If $L$ is generated by a regular sequence in $A$,
then $\pd_AG(I)$ is finite
and the equivalent conditions of Theorem \ref{2.3} hold. Hence
by part (4) of Corollary \ref{2.6}, $K$ is also
generated by a regular sequence. However, as we demonstrate in
Example \ref{2.10}, the converse fails in
general, that is, without the hypothesis that $\pd_A G(I) < \infty$,
it can happen that $K$ is generated by a regular sequence,
while  $L$ is not generated by a regular sequence.

\item Suppose  $R/I$ is Gorenstein and  $\mu(I) = \dim R + 2$. If
$G(I)$ is Gorenstein and $\pd_AG(I) < \infty$, then $L$ is
generated by an $A$-regular sequence. For $\mu(I) = \dim R + 2$
implies $\hgt(L) = 2$; then $G(I) = A/L$ Cohen-Macaulay  and
$\pd_A(A/L) < \infty$ imply $\pd_A(A/L) = 2$ by the
Auslander-Buchsbaum formula. Since $A$ and $A/L$ are Gorenstein
rings, the homogeneous minimal free resolution of $A/L$ has the
form
$$
\CD
0 @>>> A  @>>> A^2 @>>> A @>>> A/L @>>> 0.
\endCD
$$
Hence  $\mu(L)=2=\hgt(L)$, and therefore $L$
is generated by a regular sequence.

\item
Assume that $R$ is Gorenstein and
$J$ is a reduction of $I$ with $\mu(J) = \dim R$.
If $r_J(I)\leq 1$ and $\dim R \ge 2$, then it follows from
\cite[Cor.3.2]{CP} that $G(I)$ is Gorenstein if and only
if $I = J:I$. If $R/I$ is Gorenstein and $r_J(I) \le 1$,
 we prove by induction on
$\dim R$ that $G(I)$ is Gorenstein implies $\pd_A G(I) < \infty$.
Suppose $\dim R = 0$, in which case $G(I) =
R/I \oplus I$. If $G(I)$ is Gorenstein and $I \neq 0$, then $0:I=I$ by Theorem \ref{3.1}.
Thus $I$ is principal since $R$ and $R/I$ are Gorenstein,
see for instance \cite[Prop.3.3.11(b), p.114]{BH}. Therefore  $I/I^2 = I \cong
R/I$ and hence $\pd_A G(I) < \infty$. Now suppose $\dim R > 0$. With
the notation of diagram (\ref{eq1}), we may assume that $\psi(X_1) = a^*
\in G(I)$ is a regular element of $G(I)$, where $a^*$ is the
leading form of some $a \in J \backslash {\m}J$.  Therefore  $G(I)/a^*G(I)
\cong G(I/aR)$  and $\pd_A G(I) = \pd_{A/X_1A} G(I)/a^*G(I)$. We
conclude by induction that  $\pd_A G(I) < \infty$. Thus in the
case where $R/I$ is Gorenstein and $r_J(I) \le 1$, if $G(I)$ is
Gorenstein, then the equivalent conditions  of Theorem \ref{2.3}
hold.

\item
There exist
examples where $R$ is Gorenstein of dimension zero, $R/I$
is Gorenstein, $I^3 = 0$ and $G(I)$ is Gorenstein, but
$\pd_A G(I) = \infty$. To obtain an example illustrating
this consider the submaximal Pfaffians $Y^3, XZ, XY^2 + Z^2, YZ, X^2$
of the matrix
$$
\begin{bmatrix}
0  & X & 0 & 0 & Z \\
-X & 0 & Y & Z & 0 \\
0 & -Y & 0 & X & 0 \\
0 & -Z & -X & 0 & Y^2 \\
-Z & 0 & 0 & -Y^2 & 0
\end{bmatrix}.
$$
Let $k$ be a field and let
$H$ denote the ideal of the polynomial  ring
$k[X, Y, Z]$ generated by these Pfaffians. Notice that
$H$ is homogeneous with respect to the grading
that assigns   $\deg X = 0$ and $\deg Y = \deg Z = 1$.
Let $R = k[X, Y, Z]/H$, write $x$, $y$, $z$ for the images of
$X$, $Y$, $Z$ in $R$
and let $I = (y, z)R$.
Then $R$ is an Artinian Gorenstein local ring by \cite[Thm.2.1]{BE}.
Furthermore $R/I=k[X]/(X^2)$ is Gorenstein and $I^3=0$.
Finally $G(I)\cong R$ by our choice of the grading.
With $A = (R/I)[U, V]$
and $\psi : A \to G(I)\cong R$ defined by
$\psi(U) = y$ and
$\psi(V) = z$, we have
$L = \ker(\psi) = (xV, UV, V^2 + xU^2, U^3)A$.
Thus  $I/I^2$ is not free over $R/I$ and then $\pd_A G(I) = \infty$
by Theorem \ref{2.3}.
This even  provides  an example where the associated
graded ring $G(I)$ is Gorenstein and the
fiber cone $F(I) \cong k[U, V]/(UV, V^2, U^3)$ is not Gorenstein.

\item Suppose $R$ is Gorenstein. If $R/I$ is Gorenstein and of
finite projective dimension over $R$, it is shown in
\cite[Thm.2.6]{NV} that
$G(I)$ is Cohen-Macaulay implies $I$ is generated
by a regular sequence.
\end{enumerate}}
\end{remark}

\begin{exam}\label{2.10}{\em
Let $k$ be a field, $R = k[[t^3, t^4, t^5]]$ and $I =
(t^3, t^4)R$. Then $R$ is a one-dimensional
Cohen-Macaulay local domain, $\mu(I) = 2$ and
$J = t^3R$ \, is a reduction of $I$.
We have $JI \subsetneq I^2$ and $JI^2 = I^3$.
Hence $r_J(I) = 2$. If  $w$ denotes the image of
$t^5$ in $R/I$, then  $R/I \cong k[w]$, where
$w^2 = 0$. Thus $R/I$ is Gorenstein with
$\lambda(R/I) = 2$.
Consider the commutative diagram with exact rows and surjective
column maps
$$
\CD
0 @>>> L := \ker(\psi)  @>>> (R/I)[X, Y] @>\psi>> G(I)  @>>> 0 \\
@. @V\pi_1' VV        @V\pi_2' VV                @V\pi_3' VV \\
0 @>>> K := \ker(\phi) @>>> (R/\m)[X, Y] @>\phi>> F(I) @>>> 0,
\endCD
$$
where $\psi(X) = t^3 + I^2$ and $\psi(Y) = t^4 + I^2$.
It is readily checked that $wX, wY$ and $Y^3$ generate
$L = \ker(\psi)$.
In particular  $K = \ker(\phi)$ is generated by $Y^3$, so
$K$ is generated by a regular sequence, while $L$ is not generated
by a regular sequence.
Also notice that $G(I)$ and $F(I)$ both have multiplicity 3,
and  $I/I^2$ has length 2 and is not a  free module over $R/I$.}
\end{exam}

\smallskip

\section{The Gorenstein property for $G(I)$}

In this section, we establish  a necessary and sufficient condition for $G(I)$ to
be Gorenstein in terms of residuation of powers of $I$
with respect to a reduction $J$ of $I$ for which $\mu(J) = \dim R$.
We first state this in dimension zero. Among
the equivalences in Theorem \ref{3.1}, the equivalence of
(1), (3) and (5) are due to Ooishi \cite[Thm.1.5]{O}.
We include elementary direct arguments in the proof.
We use the floor function $\lfloor x \rfloor $ to denote
the largest integer which is less than or equal to $x$.

\begin{thm}\label{3.1}
Let $(R,\m)$ be a zero-dimensional Gorenstein local ring and let
$I \subseteq \m$ be an ideal of $R$. Assume that $I^r \ne 0$ and
$I^{r+1} = 0$. Let $G := G(I) = \oplus_{i=0}^r G_i$  be the
associated graded ring of $I$, and let $S := \Soc(G) =
\oplus_{i=0}^r S_i$ denote the socle of $G$. Then the following
are equivalent $:$
\begin{enumerate}
\item $G$ is a Gorenstein ring. \item $S_i = 0$ \, for $0 \leq i
\leq r-1$. \item $0:_{R}I^{r-i} = I^{i+1}$ \, for $0 \leq i \leq
r-1$. \item $0:_{R}I^{r-i} = I^{i+1}$ \, for $0 \leq i \leq
\lfloor \frac{ r-1}{2} \rfloor$. \item $\lambda(G_i) =
\lambda(G_{r-i})$ \, for $0 \leq i \leq \lfloor \frac{ r-1}{2}
\rfloor$. \item $I^r :_{R} I^{r-i} = I^i$  for $1 \leq i \leq
r-1$, and $0:_{R} I = I^r$. \item $I^{r-i}/I^r$ is a faithful
module over $R/I^i$  for $1 \leq i \leq r-1$, and $I$ is faithful
over $R/I^r$.
\end{enumerate}
\end{thm}

\demo We may assume $r > 0$. Write $k=R/\m$ and let $\mathfrak{M}$
denote the  unique maximal ideal of $G$. We first compute $S_i = 0
:_{G_i} \mathfrak{M}$  \, for $0 \leq i \leq r$. Since
$\mathfrak{M}$ is generated by $\m/I$ and $I/I^2$ it follows that
$$
S_i = (0 :_{I^i/I^{i+1}} \m/I) \cap (0 :_{I^i/I^{i+1}} I/I^2)
    = \frac{I^i}{I^{i+1}} \cap
\frac{(I^{i+1} :_{R} \m)}{I^{i+1}} \cap \frac{(I^{i+2} :_{R}
I)}{I^{i+1}}.
$$
Therefore
\begin{equation}\label{star} \tag{3.2}
S_i = \frac{I^i \cap (I^{i+1} : \m) \cap (I^{i+2} : I)}{I^{i+1}}
\quad  \text{for} \quad 0 \leq i \leq r.
\end{equation}
In particular $S_r = 0:_{I^r}\m$  because  $I^{r+1} = 0$. Note that $S_r \not = 0$.
Since $(R, \m)$ is a zero-dimensional Gorenstein local ring, we have
\begin{equation}\label{starstar} \tag{3.3}
 S_r \cong k.
\end{equation}

\smallskip

\noindent
$(1) \Longleftrightarrow (2)$ : The ring $G$ is Gorenstein if and only if $\dim_k S = 1$
if and only if $S_i = 0$ for $0 \leq i \leq r-1$, by (\ref{starstar}).

\smallskip

\noindent
$(2) \implies (3)$ : Condition (2) implies that $S = S_r \cong k$  by
(\ref{starstar}).
Hence $S = s^{*}k$ for some $0 \not = s^{*} \in S_r$.

Let $0 \leq i \leq r-1$. It is clear that
$0 : I^{r-i} \supseteq I^{i+1}$
because $I^{r+1} = 0$. To show the reverse inclusion
suppose that $0 : I^{r-i} \not\subseteq I^{i+1}$.
In this case there exists an element $z \in 0:I^{r-i}$ with
$z \in I^j \backslash I^{j+1}$ for some $0 \leq j \leq i$. Since
$0 \not = z^* \in I^j/I^{j+1}$ and $S=S_r=s^{*}k$, we can express
$s^* = z^* w^*$ for some $w^* \in I^{r-j}/I^{r-j+1}$. As $s^*
\not = 0$ it follows that $zw \not = 0$. This is impossible since
$z \in 0 : I^{r-i}$ and $w\in I^{r-j} \subseteq I^{r-i}$.

\smallskip

\noindent
(3) $\implies$ (4) : This is obvious.

\smallskip

\noindent
(4) $\implies$ (5) : For
$0 \leq i \leq \lfloor \frac{ r-1}{2} \rfloor$  we have
$$
\aligned
\lambda(G_{r-i}) &= \lambda(I^{r-i}/I^{r-i+1})\\
                 &= \lambda(R/I^{r-i+1}) - \lambda(R/I^{r-i}) \\
                 &= \lambda(0 : I^{r-(i-1)}) - \lambda(0 : I^{r-i}))
\quad \text{by \cite[Prop.3.2.12(b), p.103]{BH}} \\
                 &= \lambda(I^i) - \lambda(I^{i+1}) \quad
\text{by condition (4)}\\
                 &= \lambda(I^i/I^{i+1}) = \lambda(G_i).
\endaligned
$$

\smallskip

\noindent (5) $\implies$ (3) : If condition (5) holds for $0 \leq
i \leq \lfloor \frac{ r-1}{2} \rfloor$, then it obviously holds
for every $i$. For $0 \leq i \leq r-1$ we have
$$
\aligned
\lambda(I^{i+1}) &= \lambda(G_{i+1}) + \cdots + \lambda(G_r)\\
                 &= \lambda(G_{r-(i+1)})+ \cdots + \lambda(G_0)
\quad \text{by condition (5)} \\
                 &= \lambda(R/I^{r-i})  \\
                 &= \lambda(0:I^{r-i}). \\
 \endaligned
$$
As $I^{r+1} = 0$, we  have $I^{i+1} \subseteq 0:I^{r-i}$.
Since these two ideals have the same length, we conclude that
$I^{i+1} = 0: I^{r-i}$.

\smallskip

\noindent
$(3) \implies (6)$ : Let $1 \leq i \leq r-1$. The inclusion
$I^r : I^{r-i} \supseteq I^i$ is clear. To show $"\subseteq"$,
observe that
$$
I^r : I^{r-i} \subseteq I^{r+1} : I^{r-i+1} = 0 : I^{r-(i-1)} = I^i,
$$
where the last equality holds by condition (3).

\smallskip

\noindent $(6) \implies (2)$ : From (\ref{star})  we have  for
$0 \leq i \leq r-2$,
$$
S_i \subseteq \frac{I^{i+2} : I}{I^{i+1}} \subseteq
\frac{I^{i+2+(r-i-2)} : I^{1+(r-i-2)}}{I^{i+1}}
= \frac{I^{r} : I^{r-(i+1)}}{I^{i+1}}= \frac{I^{i+1}}{I^{i+1}}, \quad
$$
where the last equality follows from condition (6). Again by
(\ref{star}) and condition (6),
$$
S_{r-1} \subseteq \frac{0:I}{I^r} = \frac{I^r}{I^r}.
$$
Hence $S_i = 0$ for $0 \leq i \leq r-1$.

\smallskip

\noindent
$(6) \Longleftrightarrow (7)$ : This is clear.
                                         \QED

\medskip

\medskip

\setcounter{lem}{3}

\begin{cor}\label{3.3}
Let $(R,\m)$ be a zero-dimensional Gorenstein local ring
and let $I \subseteq \m$ be an ideal of $R$. Assume
that $I^r \ne 0$ and $I^{r+1} = 0$.
If   $G(I)$ is  Gorenstein, then
$I^i/I^{i+1}$ is a faithful module over $R/I$ for $0 \leq i \leq  r$.
\end{cor}

\demo
It suffices to show that $I^{i+1} : I^i = I$ for $0 \leq i \leq r$.
The inclusion $``\supseteq"$ is clear. To show the
inclusion $``\subseteq"$, observe that
we have
$$
I^{i+1} : I^i
\subseteq I^{r+1} : I^r = 0 : I^r = I,
$$
where the last equality follows from condition (3) of Theorem \ref{3.1}.
\QED

\medskip

\begin{remark} \label{3.4}
{\em
\begin{enumerate}
\item
The existence of a zero-dimensional Gorenstein local ring
$(R,\m)$ for which the associated graded ring
$G(\m)$ is not Gorenstein shows that
the converse of Corollary \ref{3.3} is not true. A specific
example illustrating this is given in Example \ref{3.5}.

\item
In general, suppose $I \subseteq \m$ is an ideal of a Noetherian
local ring $(R,\m)$. If the ideal $G(I)_+ = \oplus_{i \ge
1}(I^i/I^{i+1})$ of $G(I)$ has positive grade, then $I^i/I^{i+1}$
is a faithful module over $R/I$ for every $i \ge 0$. For if
$I^i/I^{i+1}$ is not faithful over $R/I$, then there exists $a \in
R \setminus I$ such that $aI^i\subseteq I^{i+1}$. Hence $aI^{i+j}
\subseteq I^{i+j+1}$ and $I^{i+j}/I^{i+j+1}$ is
not faithful over $R/I$ for every integer $j \geq 0$.
On the other hand, if $G(I)_+$ has
positive grade, then there exists a homogeneous $G(I)$-regular
element $b \in I^k/I^{k+1}$ for some positive integer $k$. Then
$b^\ell \in I^{k\ell}/I^{k\ell +1}$ is $G(I)$-regular, and
therefore $I^{k\ell}/I^{k\ell+1}$ is a faithful module over $R/I$
for every positive integer $\ell$.
\end{enumerate}}
\end{remark}

\medskip

\begin{exam}\label{3.5}{\em
Let $k$ be a field and $R = k[[t^5, t^6, t^9]]$. As observed by Sally
in \cite[Ex.3.6]{S},
$$
R \cong k[[X, Y, Z]]/(YZ - X^3, Z^2 - Y^3)
$$
and
$$
G(\m) \cong k[X, Y, Z]/(YZ, Z^2, Y^4 - ZX^3).
$$
Thus $(R,\m)$ is a  one-dimensional Gorenstein local domain
and $G(\m)$ is Cohen-Macaulay, but not Gorenstein.
It is readily seen that $J = t^5R$ \, is a minimal reduction of
$\m$ with $r_J(\m) = 3$. Sally also observes
that $\overline{R} = R/J$ is a zero-dimensional
Gorenstein local ring with maximal ideal
$\n = \m \overline{R}$ such that $G(\n)$ is
not Gorenstein. Indeed, the leading form $(t^5)^*$ of
$t^5$ in $\m/\m^2$ is  a regular element of $G(\m)$.
Therefore $G(\n) \cong G(\m)/(t^5)^*G(\m)$. It follows that
$G(\n)$ is not Gorenstein, $\n^3 \neq 0$ and $\n^4 = 0$.
Since $\n^i/\n^{i+1}$
is a vector space over $\overline{R}/\n$,  $\n^i/\n^{i+1}$
is a nonzero free $(\overline{R}/\n)$-module for $0 \leq i \leq 3$ and
hence in particular a faithful $(\overline{R}/\n)$-module.
The dimensions of the components of $G(\n) = R/\n \oplus
\n/\n^2 \oplus \n^2/\n^3 \oplus \n^3$ are $1, 2, 1, 1$. This
nonsymmetry reflects the fact that $G(\n)$ is not Gorenstein,
see Theorem \ref{3.1}.}
\end{exam}


We now turn to rings of arbitrary dimension. For this the following
lemma is needed.

\begin{lem}\label{3.6}
Let $(R,\m)$ be a $d$-dimensional Gorenstein local ring and let
$I$ be an $\m$-primary ideal. Assume that $J \subseteq I$ is a
reduction of $I$ with $\mu(J) = d$,  let $r_J(I)$ denote the
reduction number of $I$ with respect to $J$ and let $r$ be an
integer with $r \ge r_J(I)$. If $J:_{R} I^{r-i}  = J + I^{i+1}$
for every $i$ with $0 \le i \le r-1$, then $G(I)$ is
Cohen-Macaulay.
\end{lem}

\demo
The assertion is clear for $d = 0$. Next we assume
that  $d = 1$. By  Valabrega-Valla  \cite[Cor.2.7]{VV}, it
suffices  to verify that
$J \cap I^i = JI^{i-1}$ for $1 \le i \le r$. We first
prove that
\begin{equation}\label{claim} \tag{3.8}
I^r : I^{r-i} = I^i  \quad\quad {\rm for} \quad\quad 1 \le i \le r \,.
\end{equation}

To show (\ref{claim})
we proceed by induction on $i$. One has
$$
I^i \subseteq I^r : I^{r-i} \subseteq I^{r+1} : I^{r-i+1}
= JI^r : I^{r-(i-1)} \subseteq J : I^{r-(i-1)}
=J + I^{(i-1)+1}  = J + I^i.
$$
This proves the case $i=1$. Now assume $i \geq 2$. The containment
$I^r:I^{r-i} \subseteq J+I^i$ gives
$$
\aligned
I^r : I^{r-i} &= (J + I^i) \cap (I^r : I^{r-i}) \\
                    &= [J \cap (I^r : I^{r-i})] + I^i \quad \text{since $I^i \subseteq I^r:I^{r-i}$}\\
                    &= J[(I^r : I^{r-i}) : J] + I^i  \quad
\text{since $J$ is principal} \\
                    &= J(I^r :JI^{r-i})  + I^i \\
                    &\subseteq J(I^{r+1} : JI^{r-i+1}) + I^i \\
                    &= J(JI^r : JI^{r-(i-1)}) + I^i \\
                    &= J(I^r : I^{r-(i-1)}) + I^i  \quad
\text {since  $J$ is principal and regular} \\
                    &= JI^{i-1} + I^i \quad \quad \text{by the induction hypothesis}\\
                    &= I^i.
\endaligned
$$
This completes the proof of (\ref{claim}).

Now we have for $1 \leq i \leq r$,
$$
\aligned
J \cap I^i &= J(I^i : J) \quad \text {since $ J$ is principal} \\
           &\subseteq J(I^{r+1} : JI^{r+1-i}) \\
           &= J(JI^r : JI^{r-(i-1)}) \\
           &= J(I^r : I^{r-(i-1)}) \quad \text {since $J$ is principal and regular} \\
           &= JI^{i-1} \qquad \text {by (\ref{claim})}.
\endaligned
$$
Hence by the criterion of Valabrega-Valla, $G(I)$ is Cohen-Macaulay.
This completes the proof of Lemma \ref{3.6} in the case where $d = 1$.

Finally let $d \geq 2$. We may assume that the residue field of
$R$ is infinite. There exist elements $a_1, \ldots ,a_{d-1}$ that form
part of a minimal generating set of $J$ and a superficial sequence
for $I$. Write $\overline R = R/(a_1, \ldots, a_{d-1})R$ \, and
$\overline{I} = I \overline{R}$. Since $\dim \overline{R} = 1$,
we have that $G(\overline I)$ is a Cohen-Macaulay ring,
necessarily of dimension one. Hence by \cite[Lemma
2.2]{HM}, $G(I)$ is Cohen-Macaulay. \QED

\medskip

We are now ready to prove the main result of this section. Notice that we
do not require $G(I)$ to be Cohen-Macaulay in this theorem; instead, the
Cohen-Macaulayness is a consequence of the colon conditions (2) or (3) of
the theorem.

\setcounter{lem}{8}

\begin{thm}\label{3.8}
 Let $(R,\m)$ be a $d$-dimensional Gorenstein local
ring  and let $I$ be an $\m$-primary ideal. Assume that $J
\subseteq I$ is a reduction of $I$ with  $\mu(J) = d$, and let $r
= r_J(I)$ be the reduction number of $I$ with respect to $J$. Then
the following are equivalent $:$
\begin{enumerate}
\item $G(I)$ is Gorenstein. \item $J :_{R} I^{r-i} = J + I^{i+1}$
for $0 \leq i \leq r-1$. \item $J :_{R} I^{r-i} = J + I^{i+1}$ for
$0 \leq i \leq \lfloor \frac{ r-1}{2} \rfloor$.
\end{enumerate}
\end{thm}
\demo
The equivalence of items (2) and (3) follows from the double annihilator property
in the zero-dimensional Gorenstein local ring $R/J$, see for
instance \cite[(3.2.15), p.107]{BH}. To
prove the equivalence of (1) and (2),
by Lemma \ref{3.6}, we may assume that $G(I)$ is a Cohen-Macaulay ring.
Write $J=(a_1, \ldots, a_d)R$ and set $\overline R =R/J$,
$\overline I  =I/J$. Since $G(I)$ is Cohen-Macaulay and $J$ is
a minimal reduction of $I$, it follows that $a_1^*, \ldots, a_d^*$ form a regular
sequence on $G(I)$ and therefore
$$
G(I)/(a_1^*, \ldots, a_d^*)G(I) \cong G(I/J).
$$
In particular ${\overline I}^r \neq 0$. Hence $(\overline R, \overline \m)$
is a zero-dimensional Gorenstein local ring with
${\overline I}^r \neq 0$ and ${\overline I}^{r+1} = 0$.
Now
$$
\aligned
G(I) \quad \text {is Gorenstein}
&\Longleftrightarrow G(\overline I) \quad \text{is Gorenstein} \\
&\Longleftrightarrow 0 : {\overline I}^{r-i} = {\overline I}^{i+1}
 \quad \text {for } 0 \leq i \leq r-1  \text { by Theorem 3.1}\\
&\Longleftrightarrow J : I^{r-i} = J + I^{i+1}
\quad \text {for }  0 \leq i \leq r-1.\\
\endaligned
$$
This completes the proof of Theorem \ref{3.8}.
\QED

\medskip

We  record the following corollary to Theorem \ref{3.8}
for the case of reduction number two.

\begin{cor}\label{3.2}
Let $(R,\m)$ be a $d$-dimensional Gorenstein local ring and let
$I$ be an $\m$-primary ideal. Assume that $J \subseteq I$ is a
reduction of $I$ with $\mu(J) = d$ and that $r_J(I) = 2$, i.e.,
$JI \neq I^2$ and $JI^2 = I^3$. Then the following are equivalent
$:$
\begin{enumerate}
\item $ G(I)$ is  Gorenstein. \item $J:_{R} I^2 = I$.
\end{enumerate}
\end{cor}

\smallskip

The next corollary to Theorem \ref{3.8}  deals with the
problem of lifting the Gorenstein property of associated
graded rings.  Notice we are not assuming that $G(I)$ is
Cohen-Macaulay.

\begin{cor}\label{3.9}
Let $(R,\m)$ be a $d$-dimensional Cohen-Macaulay  local ring
and let $I$ be an $\m$-primary ideal. Assume that
$J \subseteq I$ is a reduction of $I$ with $\mu(J) = d$
and that $I^r \not\subseteq J$ for $r = r_J(I)$ the reduction
number of $I$ with respect to $J$. Set $\overline{I}: = I/J
\subseteq \overline{R}: = R/J$.
If $G(\overline{I})$ is Gorenstein, then $G(I)$ is Gorenstein.
\end{cor}
\demo If $G(\overline{I})$ is Gorenstein, then $\overline{R}$ and
hence $R$ are  Gorenstein \cite[p.121]{M}. Since $I^r
\not\subseteq J$, we have that  $r= r_J(I)$ is also the reduction
number of $\overline{I}$ with respect to the zero ideal. Now the
assertion follows from Theorem \ref{3.8}. \QED

\medskip

In Example \ref{3.11} we exhibit a one-dimensional Gorenstein local
domain $(R,\m)$,  an $\m$-primary ideal $I$ and a principal
reduction $J$ of $I$ such that  for
$\overline{I}: = I/J \subseteq \overline{R}: = R/J$, the associated
graded ring $G(\overline{I})$ is Gorenstein, while $G(I)$ is not Gorenstein.

\smallskip

\begin{exam}\label{3.11}{\em
Let $k$ be a field, $R = k[[t^4, t^5, t^6]]$ and $I =
(t^4, t^5)R$. Then $R$ is a one-dimensional
Gorenstein local domain, $\mu(I) = 2$ and
$J = t^4R$ is a principal reduction of $I$.
An easy computation shows that the reduction number $r_J(I) = 3$.
On the other hand,  $I^2 \subseteq J$. Hence the associated
graded ring $G(I)$ is not Cohen-Macaulay. We have
$\lambda(R/I) = \lambda(\overline{R}/\overline{I}) = 2$ and
$\lambda(\overline{R}) = 4$. Therefore,  by Theorem \ref{3.1},
$G(\overline{I}) = \overline{R}/\overline{I} \oplus \overline{I}$ is
Gorenstein.
On the other hand, the
principal reduction $J' = (t^4 - t^5)R$ of $I$ has the
property that $I^2$ is not contained in $J'$. Therefore the  associated
graded ring $G(I/J')$ has Hilbert function $2, 1, 1$ and
thus is not Gorenstein. }
\end{exam}

\medskip

\section{The Quasi-Gorenstein property of the extended Rees algebra }

In Theorem \ref{4.1} we give a general
characterization for when the extended Rees algebra
 $R[It, t^{-1}]$ is   quasi-Gorenstein.
In case $G(I)$ is Cohen-Macaulay, this characterization would also
follow from \cite[Thm.5.3]{GI}.  In fact, if $G(I)$ is
Cohen-Macaulay, then the quasi-Gorensteinness of $R[It, t^{-1}]$
is equivalent to the Gorensteinness of $G(I)$, since a
Cohen-Macaulay quasi-Gorenstein ring is Gorenstein and $G(I) =
R[It, t^{-1}]/(t^{-1})$  is Cohen-Macaulay (resp. Gorenstein)  if
and only if $R[It, t^{-1}]$ is Cohen-Macaulay (resp. Gorenstein).

For an ideal $I$ of a ring $R$ and an integer $i$,
we define $I^i = R$ if $i \le 0$.

\medskip

\begin{thm}\label{4.1} Let $(R,\m)$ be a Gorenstein local ring
and let $I \subseteq \m$ be an ideal with $\hgt(I) = g > 0$.
Assume that $J \subseteq I$ is a reduction of $I$ with
$\mu(J) = g$. Let $k$ be an integer with
$k \ge r:=r_J(I)$,  the reduction number
of $I$ with respect to $J$,  and let $B = R[It, t^{-1}]$ be
the extended Rees algebra of the ideal $I$. Then the graded
canonical module $\omega_B$ of $B$ has the form
$$
\omega_B = \bigoplus_{i \in \Z}(J^{i+k}:_RI^k)t^{i+g-1}.
$$
In particular, for $a \in \Z$,  the following are equivalent $:$
\begin{enumerate}
\item
$R[It, t^{-1}]$ is quasi-Gorenstein with {\bf a}-invariant $a$.
\item
$J^i:_RI^k = I^{i+a-k+g-1}$ for every $i \in \Z$.
\end{enumerate}
\end{thm}

\demo  Let $K = \Quot(R)$ denote the total ring of quotients of $R$
and let $A = R[Jt, t^{-1}] \subseteq B = R[It, t^{-1}]$.
Notice that $A/(t^{-1}) \cong G(J)$, where $t^{-1}$ is a
homogeneous $A$-regular element of degree $-1$. Moreover,
since $J$ is generated by a regular sequence, $G(J)$ is a standard
graded polynomial ring in $g$ variables over the Gorenstein
local ring $R/J$. Thus $A$ is Cohen-Macaulay and
$\omega_A \cong A(-g+1) \cong At^{g-1}$.

The extension $A \subseteq B$ is finite and
$\Quot(A) = \Quot(B) = K(t)$ since $g > 0$. Therefore
\[
\omega_B \cong \Hom_A(B,\omega_A) \cong
At^{g-1}:_{K(t)}B = (A:_{K(t)}B)t^{g-1}
= (A :_{R[t, t^{-1}]}  B)t^{g-1},
\]
where the last equality holds because
$$
A:_{K(t)}B \subseteq A:_{K(t)} 1 \subseteq A
\subseteq R[t,t^{-1}].
$$
We may now make the identification
$\omega_B = (A:_{R[t,t^{-1}]}B)t^{g-1}$.

Let $i$ and $j$ be any integers. Since $J$ is a complete
intersection, it follows that $J^{i+j+1}:_R J = J^{i+j}$.
Hence
$$
J^{i+j}:_R I^j = (J^{i+j+1}:_RJ):_RI^j = J^{i+j+1}:_RJI^j
\supseteq J^{i+j+1}:_RI^{j+1},
$$
where the last inclusion is an equality whenever $j \ge k \ge r_J(I)$. Therefore
$$
\cap_{j \in \Z}(J^{i+j}:_R I^j) = J^{i+k}:_RI^k.
$$

We conclude that

$$
[A:_{R[t,t^{-1}]}B]_i = (\cap_{j \in \Z}(J^{i+j}:_R I^j))t^i
= (J^{i+k}:_RI^k)t^i,
$$
which gives
$$
\omega_B = \bigoplus_{i \in \Z}(J^{i+k}:_RI^k)t^{i+g-1}.
$$
This description shows in particular that $[\omega_B]_i = Rt^i$ for
$i << 0$.

Now $B$ is quasi-Gorenstein with {\bf a}-invariant $a$ if and only
if $\omega_B \cong B(a)$, or equivalently, $\omega_B = \alpha Bt^{-a}$
for some unit $\alpha$ of $K$. As $[\omega_B]_i = Rt^i$ for $i << 0$,
it follows that $\alpha$ is necessarily a unit in $R$. Thus
$\omega_B \cong B(a)$ if and only if
$$
\oplus_{i \in \Z}(J^{i+k}:_R I^k)t^{i+g-1} = \oplus_{i \in \Z}I^it^{i-a},
$$
which means that $J^i:_RI^k = I^{i+a-k+g-1}$ whenever $i \in \Z$.\QED

\medskip

\begin{cor} \label{4.2}
With notation as in Theorem $\ref{4.1}$, the following are
equivalent $:$
\begin{enumerate}
\item
The extended Rees algebra $R[It, t^{-1}]$ is quasi-Gorenstein.
\item There exists an  integer $u$ such that  $J^i:_{R} I^r =
I^{i-u}$ for every $i \in \Z$.
\end{enumerate}
If these equivalent conditions hold, then $u = r -g +1 -a \geq 0$ with
$a$ denoting the $\bf a$-invariant of $R[It,t^{-1}]$.
\end{cor}

\demo  To prove the last assertion notice that $u$ is uniqely determined since
$I$ is not nilpotent. Thus  part (2) of Theorem \ref{4.1} gives $u=r-g+1-a$. The
inequality $u \geq 0$ can be seen by setting $i=0$ in part (2) of the
present corollary.  \QED

\medskip

\begin{cor} \label{4.21}
With notation as in Theorem $\ref{4.1}$, $I^n(J^i:_{R} I^r)
\subseteq J^{n+i}:_{R} I^r$ for all integers $n$, $i$, and
$J^n(J^i:_{R} I^r) = J^{n+i}:_{R} I^r$ for every $n \geq 0$ and $i
>> 0$.
\end{cor}

\demo
This is clear since
$ \omega_B = \oplus_{i \in \Z}(J^{i+r}:_RI^r)t^{i+g-1} $
is a graded module over $B = R[It, t^{-1}]$ that  is finitely
generated over $A = R[Jt, t^{-1}]$. \QED

\medskip

\begin{remark} \label{4.3} {\em
In the setting of Theorem \ref{4.1}  we define the
{\it index of nilpotency } of $I$ with respect to $J$ to be
$s_J(I)  = \min\{i \, | \, I^{i + 1} \subseteq J \}$.

Suppose $R[It, t^{-1}]$ is quasi-Gorenstein and let $u = r-g+1-a \geq 0$ be
as in Corollary \ref{4.2}.
The corollary implies $J:I^r = I^{1-u}$ and hence $I^{1-u+r} \subseteq
J$. Thus $1-u+r \ge s_J(I) + 1$, so $r_J(I) - s_J(I) \ge u \geq 0$.
Therefore the ${\bf a}$-invariant $a$ of $R[It, t^{-1}]$ satisfies
$s_J(I) - g + 1 \le a \le r_J(I) -g + 1$. Let $w = \max\{n \, | \,
I^r \subseteq J^n \}$. Then $w = u$; for $I^{w-u}=J^w:I^r = R$
implies $u \ge w$, while $I^r \not\subseteq J^{w+1}$ implies
$I^{w+1-u}= J^{w+1}:I^r \subsetneq  R$, so $w \ge u$. It follows
that  $s_J(I) =  r_J(I)$ if and only if  $u = 0$ if and only
if $a = r_J(I) - g +1$. These
equalities hold in case $G(I)$, or equivalently $R[It,t^{-1}]$,
is  Cohen-Macaulay. }
\end{remark}

\medskip

Since in dimension two, quasi-Gorenstein is equivalent to
Gorenstein, we have the following corollary:

\begin{cor} \label{4.4}
Let $(R,\m)$ be a one-dimensional Gorenstein local ring and let
$I$ be an $\m$-primary ideal. Assume that $J$ is a principal
reduction of $I$ and that $r = r_J(I)$ is the reduction number of
$I$ with respect to $J$. Then the following are equivalent $:$
\begin{enumerate}
\item $G(I)$ is Gorenstein. \item $R[It, t^{-1}]$ is Gorenstein.
\item There exists an integer $u$ such that $J^i:_{R} I^r =
I^{i-u}$ for every $i \in \Z$. \item $J^i:_{R} I^r = I^i$ \, for
every $i \in \Z$. \item $J^i:_{R} I^r = I^i$ \, for $1 \leq i \leq
r$.
\end{enumerate}
\end{cor}

\demo
It remains to prove that (5) implies (4). This
follows since for $i > r$
we have $I^i = J^{i-r}I^r = J^{i-r}(J^r:I^r)
= J^i:I^r$,
where the last equality holds because $J$ is
principal generated by a regular element.
\QED

\begin{remark} \label{4.6} {\em
With notation as in Corollary \ref{4.4}, we have:
\begin{enumerate}
\item  As is well-known in the one-dimensional case, the reduction
number $r=r_J(I)$ is independent of the principal reduction
$J=xR$. Furthermore the ideal $J^r:_RI^r$ is also independent of
$J$. Indeed, the subring $R[I/x] = I^r/x^r =
\bigcup_{i=1}^\infty(I^i:_KI^i)$ of $K : =\Quot(R)$ is the blowup
of $I$, so is independent of $J$, and $J^r:_RI^r=J^r:_KI^r$ is the
conductor of $R[I/x]$ into $R$. \item If $I^r \subseteq J$, then
$I^r \subseteq J'$ for every reduction $J'$ of $I$, that is, $I^r$
is contained in the core of $I$. Indeed, we have $r>0$ and then
$I^r \subseteq J \iff J^{r-1}I^r \subseteq J^r \iff I^{2r-1}
\subseteq J^r \iff I^{r-1} \subseteq  J^r:I^r$. Notice that the
latter ideal is independent of $J$ by item (1). Therefore  if the
index of nilpotency $s_J(I) = r$ for one principal reduction $J$,
then $s_{J'}(I) = r$ for every principal reduction $J'$ of $I$. On
the other hand, if $s_J(I) < r$ there may exist a principal
reduction $J'$ of $I$ such that $s_J(I) \ne s_{J'}(I)$ as is
illustrated in Example \ref{3.11}. \item If $J^r:I^r = I^r$, then
$I^r \not\subseteq J$ and $s_J(I) = r$. Indeed, we may assume
$r>0$,  hence $I^{r-1}  \not\subseteq I^r$ and by item (2),
$I^{r-1} \not\subseteq  J^r:I^r \iff I^r \not\subseteq J$.
\end{enumerate}
}
\end{remark}

\medskip

Let $R$ be a Noetherian ring and $I$ an $R$-ideal containing a
regular element. The ideal $\widetilde I =
\bigcup_{i=1}^\infty(I^{i+1}:_RI^i)$ first studied by Ratliff and
Rush in \cite{RR} is called the {\it Ratliff-Rush ideal associated
to}  $I$, and $I$ is said to be a {\it Ratliff-Rush ideal}  if $I
= \widetilde I$. It turns out that the ideal $G(I)_+$ of $G(I)$
has positive grade if and only if all powers of $I$ are
Ratliff-Rush ideals \cite[(1.2)]{HLS}. Thus if $\dim R = 1$, then
$G(I)$ is Cohen-Macaulay if and only if all powers of $I$ are
Ratliff-Rush ideals.

\begin{remark} \label{4.9} {\em
Let $(R,\m)$ be a Gorenstein local ring and $I \subseteq \m$ an
ideal with $\hgt(I) = g > 0$. As in Theorem \ref{4.1}, assume
there exists a reduction $J$  of $I$ with $\mu(J) = g$. Let $C =
\bigoplus_{i \in \Z}\widetilde{I^i}t^i$ denote the extended Rees
algebra of the Ratliff-Rush filtration associated to $I$, and let
$k \geq 0$ be an integer such that $J^j \widetilde{I^k}=
\widetilde{I^{j+k}}$ for every $j \geq 0$. The following are
equivalent $:$
\begin{enumerate}
\item $C$ is quasi-Gorenstein with ${\bf a}$-invariant $b$. \item
$J^i:_{R} \widetilde{I^k} = \widetilde{I^{i+b-k+g-1}}$ for every
$i \in \Z$.
\end{enumerate}
To show this equivalence one
proceeds as in the proof of Theorem \ref{4.1}
}
\end{remark}

\begin{cor} \label{4.61}
With notation as in Corollary  $\ref{4.4}$, if $J^r:_{R} I^r =
I^{r-u}$ for some integer $u \ge 0$, then $u = 0$. If, in
addition, all powers of $I$ are Ratliff-Rush ideals, then $G(I)$
is Gorenstein.
\end{cor}

\demo
The equality
$J^r:I^r = I^{r-u}$ implies $I^{r+1-u} = I(J^r:I^r)$
and by Corollary \ref{4.21}, $I(J^r:I^r)
\subseteq J^{r+1}:I^r$. On the other hand,
as $J$ is generated by a regular element,
$J^{r+1}:I^r =  J(J^r:I^r) = JI^{r-u}$.
Therefore $I^{r+1-u} = JI^{r-u}$,  and $u = 0$
since $r$ is minimal such that $I^{r+1} = JI^r$.
By Corollary \ref{4.4}, to show $G(I)$ is Gorenstein
it suffices to show
$J^i:I^r = I^i$ for $1 \leq i \leq r$.
Again according to Corollary \ref{4.21}, $I^{r-i}(J^i:I^r) \subseteq J^r:I^r$.
Since $J^r:I^r = I^r$, it follows that
$J^i:I^r  \subseteq I^r:I^{r-i}$. We always have
$I^r:I^{r-i} \subseteq \widetilde{I^i}$. Therefore,  if in
addition,  $\widetilde{I^i} = I^i$, then  $G(I)$ is  Gorenstein.
\QED

\medskip

With notation as in Corollary \ref{4.4}, it can happen that
$J^r:_{R} I^r = I^r$ and yet $G(I)$ is  not Cohen-Macaulay.  We
illustrate this in Example \ref{4.5}.

\begin{exam}\label{4.5}{\em
Let $k$ be a field, $R = k[[t^5, t^6, t^7, t^8]]$ and $I =
(t^5, t^6, t^7)R$. Then $R$ is a one-dimensional
Gorenstein local domain with integral closure
$\overline R = k[[t]]$ and
$J = t^5R$ is a reduction of $I$.
We have $JI \subsetneq I^2 = \m^2 = t^{10}\overline R$ and $JI^2 = I^3$.
Hence $r_J(I) = 2$. Since $t^8 \not\in I$, but
$t^8I \subseteq I^2$, $G(I)$ is not Cohen-Macaulay.
To see that $J^2:I^2 = I^2$, observe that
$t^{19} \not\in J^2$ and that for each integer $i$
with $5 \le i \le 8$ we have
$t^{19-i} \in I^2$ and
$t^{19} = t^{19-i}t^i$. Notice also that the
Ratliff-Rush ideal $\widetilde{I}$ associated to
$I$ is $\m = (t^5, t^6, t^7, t^8)R$, $r_J(\m) = r_J(I) = 2$,
$J^i:\m^2 = \m^i$ for every $i$ (equivalently, for
$1 \leq i \leq 2$), and hence $G(\m)$ is
Gorenstein by Corollary \ref{4.4}.
}
\end{exam}

\begin{exam}\label{4.51}{\em
Let $R = k[[t^4, t^5, t^6]]$,  $I = (t^4, t^5)R$ \,
and $J = t^4R$ \, be as
in Example \ref{3.11}.
The Ratliff-Rush ideal
$\widetilde{I}$ associated to $I$ is $\m=(t^4,t^5,t^6)R$,
and in fact
$\widetilde{I^i} = \m^i$ for every $i \in \Z$.
We have $r_J(I) = 3$ while $r_J(\m) = 2$.  Also
$I^3 = \m^3$. We have
$J^i:I^3 = \m^{i-1}$ while $J^i:\m^2 = \m^i$
for every $i \in \Z$. In particular $G(\m)$
is Gorenstein.
}
\end{exam}

\begin{ques}\label{4.7} {\em
With notation as in  Theorem \ref{4.1}, is
the extended Rees algebra $R[It, t^{-1}]$ Gorenstein
if it is quasi-Gorenstein?  Equivalently, is the
associated graded ring $G(I)$ then Gorenstein? }
\end{ques}

If $\dim R = 2$ and $R$ is pseudo-rational in the sense of
\cite[p.102]{LT}, we observe an affirmative answer to
Question~\ref{4.7}.

\begin{cor} \label{4.8}
Let $(R, \m)$  be a $2$-dimensional pseudo-rational Gorenstein
local ring and let $I$ be an $\m$-primary ideal.
If $B = R[It, t^{-1}]$ is quasi-Gorenstein,
then $B$ is Gorenstein. In particular, if $R$
is a $2$-dimensional regular local ring, then every
extended Rees algebra $R[It, t^{-1}]$ that is
quasi-Gorenstein is Gorenstein.
\end{cor}

\demo We may assume that $R/\m$ is infinite. Let
$J \subseteq I$ be a minimal reduction of $I$
and let $r = r_J(I)$. If $J = I$, then $I$ is
generated by a regular sequence and $R[It, t^{-1}]$ is
Gorenstein. Thus we may assume $J \subsetneq I$.
By \cite[Cor.5.4]{LT}, \,
$\overline{I^{r+1}} \subseteq J^r$, where
$\overline{I^{r+1}}$ denotes the integral closure
of $I^{r+1}$. In particular $\overline{I}I^r \subseteq
J^r$, which gives $\overline{I} \subseteq J^r:I^r$.
By Corollary \ref{4.2}, there exists
an integer $u$ such that $J^i:I^r = I^{i-u}$
for every $i \in \Z$. We have $u < r$ according to Remark \ref{4.3}
since $J \neq I$. Thus $J^r:I^r \subseteq J^{u+1}:I^r$.
Therefore $\overline{I} \subseteq J^r:I^r \subseteq J^{u+1}:I^r = I$,
which shows that $I = \overline{I}$ is integrally closed.
Now \cite[Cor.5.4]{LT} implies $r \le 1$. Therefore
$B$ is Cohen-Macaulay by \cite[Prop.3.1]{VV} and hence
Gorenstein. \QED


\bigskip
\bigskip


\begin{thebibliography}{GGP0}

\bibitem[BH]{BH}{W. Bruns and J. Herzog, {\em Cohen--Macaulay Rings},
revised edition, Cambridge University Press, Cambridge, 1998.}

\bibitem[BE]{BE}{D. Buchsbaum and D. Eisenbud, Algebra structures
for finite free resolutions, and some structure theorems for ideals
of codimension 3, Amer. J. Math. {\bf 99} (1977), 447-485.}

\bibitem[CP]{CP}{A. Corso and C. Polini, Links of prime ideals and
their Rees algebras, J. Algebra {\bf 178} (1995), 224-238.}

\bibitem[DV]{DV}{C. D'Cruz and J. Verma, Hilbert series of fiber
cones of ideals with almost minimal mixed multiplicity, J. Algebra
{\bf 251} (2002), 98-109. }

\bibitem[F]{F}{D. Ferrand, Suite r\'eguli\`ere et intersection
compl\`ete, C. R. Acad. Sci. Paris {\bf 264} (1967), 427-428.}

\bibitem[GI]{GI}{ S. Goto and S. Iai, Embeddings of certain graded
rings into their canonical modules, J. Algebra {\bf 228} (2000), 377-396.}

\bibitem[HLS]{HLS}{W. Heinzer, D. Lantz and K. Shah, The Ratliff-Rush ideals
in a Noetherian ring, Comm. Algebra {\bf 20} (1992), 591-622.}


\bibitem[HSV]{HSV}{J. Herzog, A. Simis and W. Vasconcelos,
Approximation complexes and blowing-up rings, J. Algebra
{\bf 74} (1982), 466-493.}

\bibitem[HH]{HH}{S. Huckaba and C. Huneke, Powers of ideals having
small analytic deviation, Amer. J. Math. {\bf 114} (1992), 367-403.}

\bibitem[HM]{HM}{S. Huckaba and T. Marley, Hilbert coefficients and the depths
of associated graded rings, J. London Math. Soc. {\bf 56} (1997), 64-76.}

\bibitem[K]{K}{E. Kunz,
{\em Introduction to Commutative Algebra and Algebraic Geometry},
Birkh\"{a}user,  Bosten, 1985.}

\bibitem[LT]{LT}{J. Lipman and B. Teissier, Pseudo-rational local
rings and a theorem of Brian\c{c}on-Skoda about integral closures
of ideals, Michigan Math. J. {\bf 28} (1981), 97-116.}

\bibitem[M]{M}{H. Matsumura, {\em Commutative ring theory},
Cambridge Univ. Press, Cambridge, 1986.}


\bibitem[NV]{NV}{S. Noh and W. Vasconcelos, The $S_2$-closure of a
Rees algebra, Results Math. {\bf 23} (1993), 149-162.}

\bibitem[NR]{NR}{D.G. Northcott and D.  Rees, Reductions of ideals
in  local rings,
Proc. Camb. Phil. Soc. {\bf 50} (1954), 145-158.}

\bibitem[O]{O}{A. Ooishi, On the Gorenstein property of the
associated graded ring and the Rees algebra of an ideal,
J. Algebra {\bf 155} (1993), 397-414.}


\bibitem[PV]{PV}{F. Planas-Vilanova, On the module of effective
relations of a standard algebra, Math. Proc. Camb. Phil. Soc. {\bf 124}
(1998), 215-229.}

\bibitem[RR]{RR}{L. J. Ratliff and D. E. Rush, Two notes on reductions
of ideals, Indiana Math. J. {\bf 27} (1978), 929-934.}

\bibitem[S]{S}{J. D. Sally, Tangent cones at Gorenstein singularities,
Compositio Math. {\bf 40} (1980), 167-175.}

\bibitem[Sh]{Sh}{K. Shah, On the Cohen-Macaulayness of the fiber cone
of an ideal, J. Algebra {\bf 143} (1991), 156-172.}

\bibitem[VV]{VV}{P. Valabrega and G. Valla, Form rings and regular sequences,
Nagoya Math. J. {\bf 72} (1978), 93-101.}

\bibitem[Va]{Va}{G. Valla, On the symmetric and Rees algebras of an ideal,
Manuscripta Math. {\bf 30} (1980), 239-255.}

\bibitem[V]{V}{W. Vasconcelos, Ideals generated by $R$-sequences,
J. Algebra {\bf 6} (1967), 309-316.}



\end{thebibliography}
\end{document}